\documentclass{amsart}
\usepackage{fullpage}
\usepackage[pdftex]{graphicx}

\usepackage[all,cmtip]{xy} 

\usepackage{amssymb,amsfonts,amsmath}

\newcommand{\iinfty}{{\mathchoice
{\begin{minipage}{.15in}\includegraphics[width=.15in]{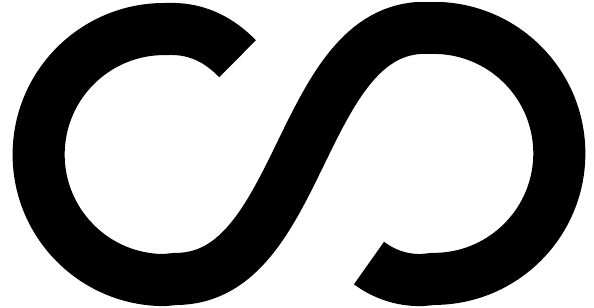}\end{minipage}}
{\begin{minipage}{.10in}\includegraphics[width=.10in]{infty2.pdf}\end{minipage}}
{\begin{minipage}{.08in}\includegraphics[width=.08in]{infty2.pdf}\end{minipage}}
{\begin{minipage}{.08in}\includegraphics[width=.08in]{infty2.pdf}\end{minipage}}
}}

\newtheorem{thm} {Theorem} 
\newtheorem{cor}[thm] {Corollary} 
 
\newtheorem{conj}[thm] {Conjecture} 
\newtheorem{prop}[thm] {Proposition} 
\newtheorem{defn}[thm]  {Definition}

\renewcommand{\int}{\operatorname{int}} 
 
\newcommand{\Ker}{\operatorname{Ker}} 
\newcommand{\Cok}{\operatorname{Cok}} 
\newcommand{\Arf}{\mathrm{Arf}}

\newcommand\W{\text{\sf W}}

\newcommand\sD{\text{\sf D}}
\newcommand\sL{\text{\sf L}}

\newcommand{\Z}{\mathbb{Z}} 
\newcommand{\N}{\mathbb{N}} 
\newcommand{\R}{\mathbb{R}}

\newcommand{\Q}{\mathbb{Q}}

\newcommand{\bW}{\mathbb{W}} 
\newcommand{\bL}{\mathbb{L}}

\newcommand{\cT}{\mathcal{T}} 
 
\newcommand{\cW}{\mathcal{W}}

\newcommand{\sra}{\twoheadrightarrow}

\title{Higher-Order Intersections in \\ Low-Dimensional Topology}

\author{Jim Conant, Rob Schneiderman and Peter Teichner}

\begin{document}
\maketitle

\begin{abstract} We show how to measure the failure of the Whitney trick in dimension 4 by constructing higher-order intersection invariants of \emph{Whitney towers} built from iterated Whitney disks on immersed surfaces in 4--manifolds. For Whitney towers on immersed disks in the 4--ball, we identify some of these new invariants with previously known link invariants like Milnor, Sato-Levine and Arf invariants. We also define higher-order Sato-Levine and Arf invariants and show that these invariants detect the obstructions to framing a \emph{twisted} Whitney tower. Together with Milnor invariants, these higher-order invariants are shown to classify the existence of (twisted) Whitney towers of increasing order in the 4--ball.
A conjecture regarding the non-triviality of the higher-order Arf invariants is formulated, and related implications for filtrations of string links and 3-dimensional homology cylinders are described.
\end{abstract}

%\title

%% When adding keywords, separate each term with a straight line: |
%\keywords{Whitney tower | Higher-order intersection | 4--manifold | link concordance | trivalent tree | Arf invariant | k-slice | Gropes}

%% Optional for entering abbreviations, separate the abbreviation from
%% its definition with a comma, separate each pair with a semicolon:
%% for example:
%% \abbreviations{SAM, self-assembled monolayer; OTS,
%% octadecyltrichlorosilane}

% \abbreviations{}

%% The first letter of the article should be drop cap: \dropcap{}
%\dropcap{I}n this article we study the evolution of ''almost-sharp'' fronts

%% Enter the text of your article beginning here and ending before
%% \begin{acknowledgements}
%% Section head commands for your reference:
%% \section{}
%% \subsection{}
%% \subsubsection{}
\begin{figure}[ht!]
        \centerline{\includegraphics[width=.5\columnwidth]{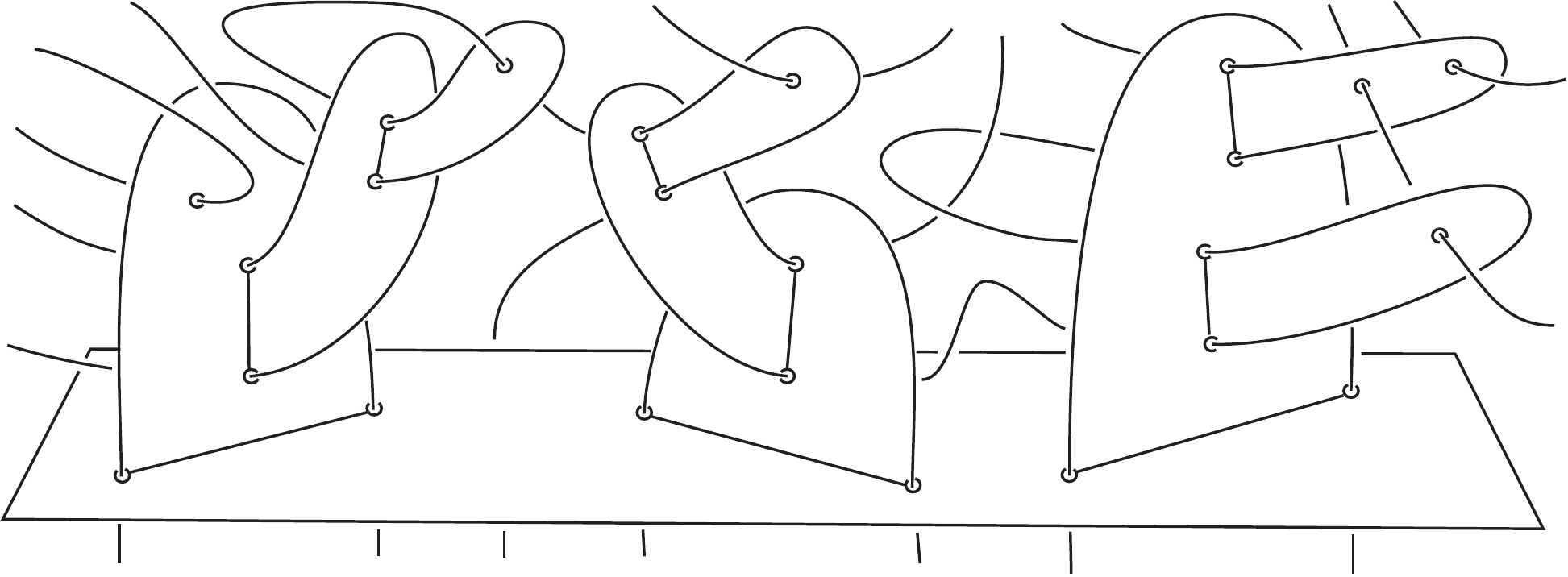}}
        \caption{Part of a Whitney tower in $4$--space.}
        \label{fig:W-tower5-fig}
\end{figure}

\section{Introduction}
Despite how it may appear in high school, mathematics is not all about manipulating numbers or functions in more and more complicated algebraic or analytic ways. In fact, one of the most interesting quests in mathematics is to find a good notion of {\em space}. It should  be general enough to cover many real life situations and at the same time sufficiently specialized so that one can still prove interesting properties about it. A first candidate was Euclidean $n$-space $\R^n$, consisting of $n$-tuples of real numbers. This covers all dimensions $n$ but is too special: the surface of the earth, mathematically modelled by the $2$--sphere $S^2$, is 2-dimensional but compact, so it can't be $\R^2$. However, $S^2$ is {\em locally} Euclidean: around every point one can find a neighborhood which can be completely described by two real coordinates (but {\em global} coordinates don't exist). 

This observation was made into the definition of an {\em $n$-dimensional manifold} in 1926 by Kneser: It's a (second countable) Hausdorff space which looks locally like $\R^n$. The development of this definition started at least with Riemann in 1854 and important contributions were made by Poincar\'e and Hausdorff at the turn of the 19th century.  It is believed to cover many important important physical notions, like the surface of the earth, the universe, and space-time (for $n=2,3,$ and $4$, respectively) but is special enough to allow interesting structure theorems. One such statement is {\em Whitney's (strong) embedding theorem}: 
Any $n$-manifold $M^n$ can be embedded into $\R^{2n}$ (for all $n\geq 1$). 

%%%%%%%%%%%%%%%%%%%%%%%%%%%%%%%%%%%%%%%%%%%
\begin{figure}[ht!]
        \centerline{\includegraphics[width=.5\columnwidth]{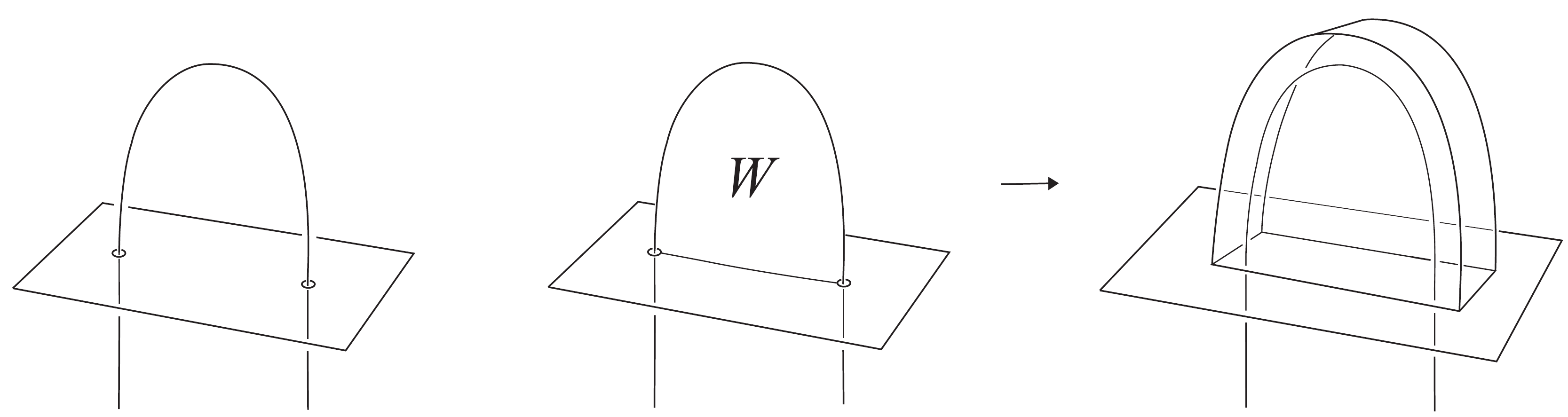}}
        \caption{Left:~A canceling pair of transverse intersections between two local sheets of surfaces in a `$3$-dimensional slice' of $4$-space. The horizontal sheet appears entirely in the `present', and the other sheet appears as an arc which is assumed to extend into `past' and `future'.
         Middle:~A Whitney disk $W$ pairing the intersections. Right:~A Whitney move guided by $W$ eliminates the intersections.}
        \label{fig:canceling-pair-and-whitney-disk-and-whitney-move}
\end{figure}
%%%%%%%%%%%%%%%%%%%%%%%%%%%%%%%%%%%%%%%%%%%%

The proof in small dimensions $n=1,2$ is fairly elementary and special, but in all dimensions $n>2$, Hassler Whitney \cite{Wh} found the following beautiful argument:  By general position, one finds an immersion $M \rightarrow \R^{2n}$ with at worst transverse double points. By adding local cusps, one can assume that all double points can be paired up by {\em Whitney disks} as in Figure~\ref{fig:canceling-pair-and-whitney-disk-and-whitney-move}, using the fact that $\R^{2n}$ is simply connected. Since $2+2<2n$ and $n+2 < 2n$, one can arrange that all Whitney disks are disjointly embedded and also meet the image of $M$ only on the boundary (as well as satisfying a certain normal framing condition). Then a sequence of (what today are called) {\em Whitney moves} leads to the desired embedding of $M$ (Figure~\ref{fig:canceling-pair-and-whitney-disk-and-whitney-move}). To be more precise, one needs to distinguish between {\em topological} and {\em smooth} manifolds. A topological $n$--manifold is locally {\em homeomorphic} to $\R^n$, whereas a smooth manifold is locally \emph{diffeomorphic} to it (in the given smooth structure). Whitney's argued in the smooth setting where transversality and isotopy extension theorems were readily available. Kirby and Siebenmann \cite{KS} made these tools available also in the topological category in dimensions $>4$, see for example page 122 for the topological Whitney move.

The Whitney move, sometimes also called the {\em Whitney trick}, remains a primary tool for turning algebraic information (counting double points) into geometric information (existence of embeddings). It was successfully used in the classification of manifolds of dimension $>4$, specifically in Smale's celebrated h-cobordism theorem \cite{Sm} (implying the Poincar\'e conjecture) and the surgery theory of Kervaire-Milnor-Browder-Novikov-Wall \cite{Wa}. The failure of the Whitney move in dimension 4 is the main reason that, even today, there is no classification of 4-dimensional manifolds in sight. 

Casson realized that in the setting of the 4-dimensional h-cobordism theorem, even though Whitney disks can't always be embedded (because $2+2=4$), they always fit into what is now called a {\em Casson tower}. This is an iterated construction that works in simply connected $4$--manifolds, where one adds more and more layers of disks onto the singularities of a given (immersed) Whitney disk \cite{Ca}. In an amazing tour de force, Freedman \cite{F1,F2} showed that there is always a {\em topologically} embedded disk in a neighborhood of certain Casson towers (originally, one needed 7 layers of disks, later this was reduced to 3). This result implied the topological h-cobordism theorem (and hence the topological Poincar\'e conjecture) in dimension $4$. At the same time, Donaldson used Gauge theory to show that the smooth h-cobordism theorem is wrong \cite{Do}, and both results were awarded with a Fields medal in 1982. Surprisingly, the smooth Poincar\'e conjecture is still open in dimension $4$ -- the only remaining unresolved case.

In the non-simply connected case, even the topological classification of $4$--manifolds is far from being understood because Casson towers cannot always be constructed. See \cite{FQ,FT1,KQ} for a precise formulation of the problem and a solution for fundamental groups of subexponential growth. However, there is a simpler construction, called a {\em Whitney tower}, which can be performed in many more instances. The current authors have developed an obstruction theory for such Whitney towers in a sequence of papers \cite{CST1,CST2,CST3,CST4,CST5,CST,S,ST,ST1}. Even though the existence of a Whitney tower does not lead to an embedded (topological) disk, it is still a necessary condition. Hence our obstruction theory provides {\em higher-order (intersection) invariants} for the existence of embedded disks, spheres, or surfaces in $4$--manifolds. 

%%%%%%%%%%%%%%%%%%%%%%%%%%%%%%%%%%%%%%%%%%%
\begin{figure}[ht!]
         \centerline{\includegraphics[width=.5\columnwidth]{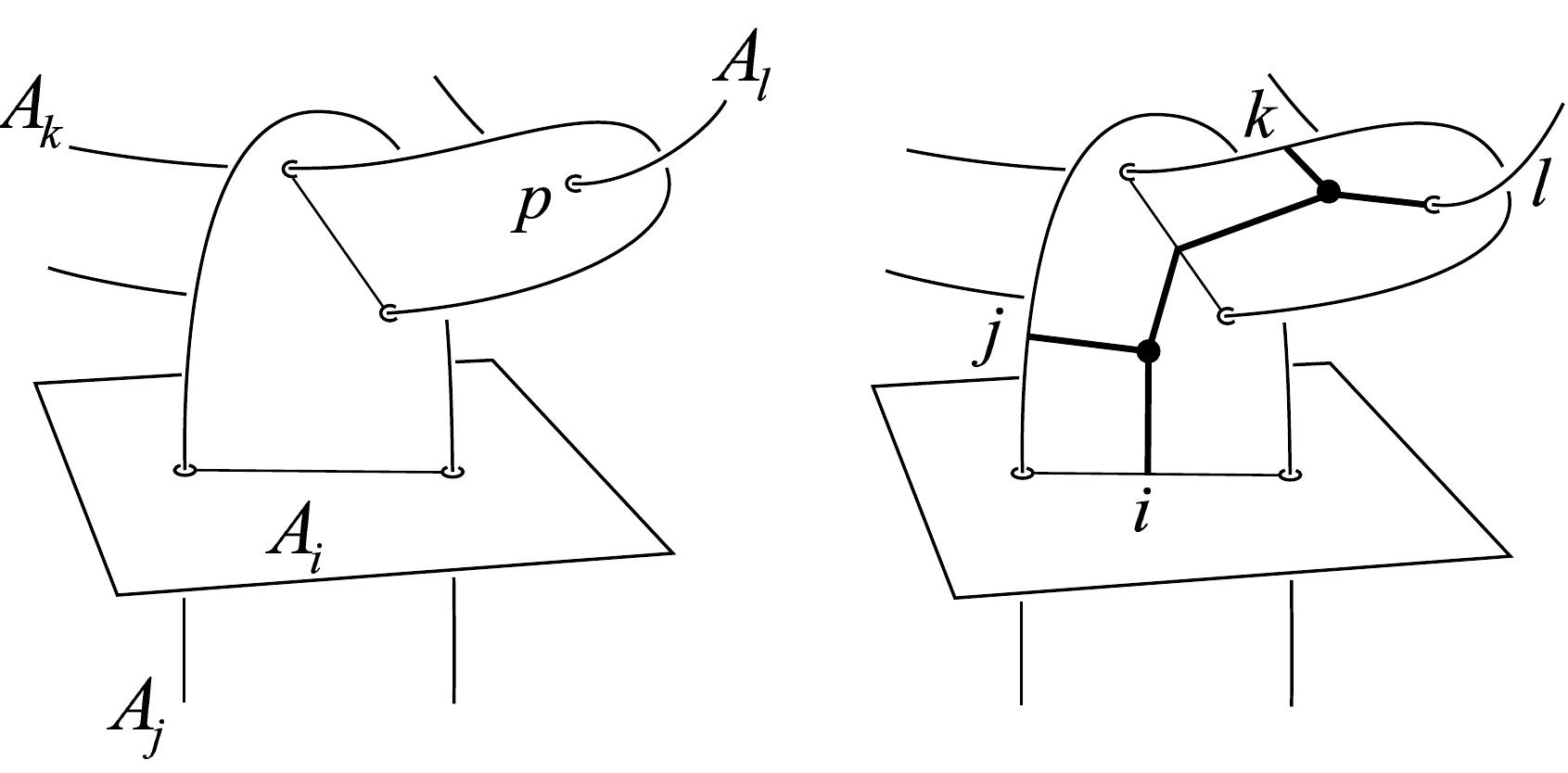}}
         \caption{Left:~Part of an order~2 Whitney tower on order~0 surfaces $A_i$, $A_j$, $A_k$, and $A_l$.
         Right:~The labeled tree $t_p$ associated to the
         order~2
         intersection point $p$.}
         \label{IHX-W-tower-fig}
\end{figure}
%%%%%%%%%%%%%%%%%%%%%%%%%%%%%%%%%%%%%%%%%%%%

The easiest example of our intersection invariant is Wall's self-intersection number for disks in 4-manifolds. If $A:(D^2,\partial D^2) \to (M^4,\partial M)$ has trivial self-intersection number  (we say that the {\em order zero invariant} $\tau_0(A)$ vanishes) then all self-intersections can be paired up by Whitney disks $W_i$. However, the $W_i$ now self-intersect and intersect each other and also the original disk $A$. Our (first order) intersection invariant $\tau_1(A,W_i)$ measures the intersections $A\pitchfork W_i$ and vanishes if they all can be paired up by (second order) Whitney disks $W_{i,j}$. This procedure continues with an invariant $\tau_2(A, W_i,W_{i,j})$ which measures both $A\pitchfork W_{i,j}$ and $W_i\pitchfork W_k$ intersections, and the construction of a higher-order Whitney tower $\cW$ if the invariant vanishes. $\cW$ is the union of $A$ (at order 0) and all Whitney disks $W_i$ (order 1), $W_{i,j}$ (order 2) and continuing with higher-order Whitney disks. If $A$ is homotopic (rel.\ boundary) to an embedding then these constructions can be continued ad infinitum.

The intersection invariants $\tau_n(A, W_i,W_{i,j},\dots)=\tau_n(\cW)$ take values in a finitely generated abelian group $\cT_n$ which is generated by certain trivalent trees that describe the 1-skeleton of a Whitney tower (Figure~\ref{IHX-W-tower-fig}).  The relations in $\cT_n$ correspond to Whitney moves, and quite surprisingly most of these relations can be expressed in terms of the so called IHX-relation which is a geometric incarnation of the Jacobi identity for Lie algebras. All the relations can be realized by controlled manipulations of Whitney towers, and as a result we recover the following approximation of the ``algebra implies geometry'' principle that is available in high dimensions:

\begin {thm}[Raising the order of a Whitney tower]\label{thm:intro-order-raising}
If $A$ supports an order $n$ Whitney tower $\cW$ with vanishing $\tau_n(\cW)$, then $A$ is homotopic (rel. boundary) to $A'$ which supports an order $n+1$ Whitney tower. Compare Theorem~\ref{thm:framed-order-raising-on-A}.
\end {thm} 

As usual in an obstruction theory, the dependence on the lower order Whitney towers makes it hard to derive explicit invariants that prevent the original disk $A$ from being homotopic to an embedding. In this paper we discuss how to solve this problem in the easiest possible ambient manifold $M=B^4$, the 4-dimensional ball. We start with maps 
\[
A_1,\dots,A_m:(D^2,S^1) \to (B^4,S^3)
\]
 which exhibit a fixed link in the boundary $3$--sphere $S^3$. If this link is \emph{slice} then the $A_i$ are homotopic (rel. boundary) to disjoint embeddings and our Whitney tower theory gives obstructions to this situation.
In the simplest example discussed above we have $m=1$ and the boundary of $A$ is just a knot $K$ in $S^3$:

\begin {thm}[The easiest case of knots \cite{CST2}]\label{thm:Arf}
The first order intersection invariant $\tau_1(A,W_i) \in \cT_1 \cong \Z_2$ can be identified with the Arf invariant of the knot $K$. It is thus a well-defined invariant that only depends on $\partial A =K$. Moreover,  it is the complete obstruction to finding a Whitney tower of arbitrarily high order $\geq 1$ with boundary $K$. 
\end {thm}

There is a very interesting refinement of the theory for knots in the setting of Cochran, Orr, and Teichner's \emph{$n$-solvable filtration}: Certain special symmetric Whitney towers of orders which are powers of $2$ have a refined measure of complexity called \emph{height}, and are obstructed by higher-order signatures of associated covering spaces \cite{COT1}. However there are no known algebraic criteria for ``raising the height'' of a Whitney tower, as given by Theorem~\ref{thm:intro-order-raising}. 

If $m>1$ then the order zero invariant $\tau_0(A_1,\dots,A_m)$ is given by the linking numbers of the components $L_i := \partial A_i$ of the link $L=\cup_{i=1}^mL_i\subset S^3$ that is the boundary of the given disks. Milnor \cite{M1,M2} showed in 1954 how to generalize linking numbers $\mu(i,j)$ inductively to higher order. Here we use the {\em total order $n$ Milnor invariants} $\mu_n$ which correspond to all length $(n+2)$ Milnor numbers $\mu(i_1,\dots,i_{n+2})$.

\begin {thm}[Milnor numbers as intersection invariants] \label{thm:pre-Milnor}
If a link $L$ bounds a Whitney tower $\cW$ of order $n$ then the Milnor invariants $\mu_k$ of order 
$k<n$ vanish. Moreover, the order $n$ Milnor invariants of $L$ can be computed from the {\em intersection invariant} $\tau_n(\cW)\in \cT_n$. Compare Theorem~\ref{thm:mu-equals-eta-of-tau}.
\end {thm}

In the remaining sections, we will make these statements precise and explain how to get {\em complete} obstructions for the existence of Whitney towers for links. Unlike the case of knots, these get more and more interesting for increasing order. In addition to the above Milnor invariants (higher-order linking numbers), we'll need higher-order versions of Sato-Levine and Arf invariants. In a fixed order, these are finitely many $\Z_2$-valued invariants, so that, surprisingly, the Milnor invariants already detect the problem up to this $2$-torsion information. 

\begin {thm}[Classification of Whitney tower concordance]\label{thm:main}
A link $L$ bounds a Whitney tower $\cW$ of order $n$ if and only if its Milnor invariants, Sato-Levine invariants and Arf invariants vanish up to order $n$. Compare Corollary~\ref{cor:mu-sl-arf-classify}.
\end {thm}

%We'll define the order $n$ Sato-Levine invariants as mod~2 reductions of certain Milnor-invariants of order $(n+1)$. To avoid this artificial shift, we also introduce the notion of a {\em twisted} Whitney tower, where certain Whitney disks are allowed to be twisted (rather than framed as is required, but not mentioned, above). Then the Milnor invariants can still be read off from an intersection invariant of a twisted Whitney tower,  and the classification is clarified as follows:
%
%
%\begin {thm}\label{thm:main-twisted}
%For $n\equiv 0,1,3 \mod 4$, a link $L$ that bounds an order $n$ twisted Whitney tower also bounds a twisted Whitney tower $\cW$ of order $n+1$ if and only if its order 
%$n$ Milnor invariants vanish. For $n\equiv 2 \mod 4$, the order $n$ Milnor and Arf invariants are the complete obstruction to the existence of a twisted Whitney tower of order $n+1$ bounded by $L$.
%\end {thm}
To prove this classification, we use Theorem~\ref{thm:intro-order-raising} to show that the intersection invariant $\tau_n(\cW)$ leads to a surjective {\em realization map} $R_n:\cT_n\twoheadrightarrow\W_n$, where $\W_n$ consist of links bounding Whitney towers of order $n$, up to order $n+1$ Whitney tower concordance (see the next section). The Milnor invariant can be translated into a homomorphism $\mu_n:\W_n\to \sD_n$, where the latter is a group defined from a free Lie algebra (which can be expressed via {\em rooted} trivalent trees modulo the Jacobi identity). The composition
\[
\eta_n: \cT_n \to \W_n \to \sD_n
\]
is hence a map between purely combinatorial objects both given in terms of trivalent trees. %(modulo the same relations).
 Using a geometric argument (grope duality), we show that it is simply given by summing over all choices of a root in a given tree (which is a more precise statement of Theorem~\ref{thm:pre-Milnor}). This map was previously studied by Jerry Levine in his work on $3$-dimensional homology cylinders \cite{L1,L2}, where he made a precise conjecture about the kernel and cokernel of $\eta_n$. He verified the conjecture for the cokernel in \cite{L3}, using a generalized Hall algorithm.

In \cite{CST3} we prove Levine's full conjecture via an application of combinatorial Morse theory to tree homology. In particular, we show that the kernel of $\eta_n$ consists only of $2$-torsion. 
This $2$-torsion corresponds to our higher-order Sato-Levine and Arf invariants, 
and is characterized geometrically in terms of a framing obstruction for \emph{twisted} Whitney towers (in which certain Whitney disks are not required to be framed). 

In the above classification of Whitney tower concordance there remains one key geometric question: Although our higher-order Arf invariants are well-defined, it is not currently known if they are in fact non-trivial. All potential values are indeed realized by simple links, so the question here is whether or not there are any further geometric relations; see Definition~\ref{def:higher-order-Arf}. We conjecture that there are indeed new higher-order Arf invariants, or equivalently, that our realization maps $\widetilde R_n: \widetilde\cT_n \to \W_n$ are isomorphisms for all $n$. Here $\widetilde \cT_n$ is a certain quotient of $\cT_n$ by what we call {\em framing relations} which come from IHX-relations on twisted Whitney towers. For 
$n\equiv 0,2,3\!\mod 4$ we do show that $\widetilde R_n$ is an isomorphism, implying that
in this further quotient the intersection invariant $\tau_n(\cW)$ only depends on the link 
$\partial \cW$, and not on the choice of Whitney tower $\cW$. The higher-order Arf invariants appear when $n=4k-3$, and our conjecture says that the same conclusion holds in these orders. 

This conjecture is in turn equivalent to the vanishing of the intersection invariants on all immersed $2$--spheres in $S^4$. Of course all such maps are null-homotopic, and a general goal of the Whitney tower theory is to extract higher-order invariants of representatives of classes in the second homotopy group $\pi_2M$. This obstruction theory is still being developed but certain aspects of it appeared in \cite{CST,S3,ST,ST1}. The fundamental group $\pi_1M$ leads to more interesting obstruction groups $\cT_n(\pi_1M)$ and a non-trivial $\pi_2M$ leads to more relations to make the intersection invariants only dependent on the order zero surfaces. 

In this paper, we will give a survey of the necessary material that is needed to understand the above results for Whitney towers in the 4-ball. More details and proofs can be found in our recent series of five papers \cite{CST1,CST2,CST3,CST4,CST5} from which we'll also survey here the following aspects of the theory:
\begin{itemize}
\item Twisted Whitney towers and obstruction to framing them
\item Geometrically $k$-slice links and vanishing Milnor invariants
\item String links and the Artin representation
\item Levine's Conjecture and filtrations of homology cylinders
\end{itemize}

%%%%%%%%%%%%%%%%%%%%%%%%%%%%%%%%%%%%%%%%%

%%%%%%%%%%%%%%%%%%%%%%%%%%%%%%%%%%%%%%%%%

\section{Whitney towers}
We work in the smooth oriented category (with discussions of orientations mostly suppressed), even though all results hold in the locally flat topological category by the basic results on topological immersions in Freedman--Quinn \cite{FQ}. In particular, our techniques do not distinguish smooth from locally flat surfaces.

Order $n$ Whitney towers are defined recursively as follows.
\begin  {defn}\label{def:w-tower}
 A \emph{surface of order 0} in an oriented $4$--manifold $M$
is a connected oriented surface in $M$ with boundary embedded in the boundary and interior immersed in the interior of $M$. 
A \emph{Whitney tower of order 0} is a collection of order 0
surfaces. The \emph{order of a (transverse) intersection point} between a surface of order $n$ and a
surface of order $m$ is $n+m$.
The \emph{order of a Whitney disk} is $(n+1)$ if it pairs intersection points of order $n$.
For $n\geq 1$, a \emph{Whitney tower of order $n$}  is a Whitney tower $\mathcal W$ of
order $(n-1)$ together with (immersed) Whitney disks pairing all order
$(n-1)$ intersection points of $\mathcal W$.
\end  {defn}

The Whitney disks in a Whitney tower may self-intersect and intersect each other as well as lower order surfaces but the boundaries of all Whitney disks are required to be disjointly embedded. In addition, all Whitney disks are required to be \emph{framed}, as will be discussed below.

\subsection{Whitney tower concordance}\label{subsec:filtr-link-conc}

We now specialize to the case $M=B^4$, and also assume that a Whitney tower $\mathcal W$ has disks for its order $0$ surfaces which have an $m$-component link in $S^3=\partial B^4$ as their boundary, denoted $\partial\mathcal W$. Let $\mathbb W_n$ be the set of all framed links $\partial \mathcal W$ where $\mathcal W$ is an order $n$ Whitney tower, and the link framing is induced by the order $0$ disks in $\cW$. This defines a filtration $\cdots\subseteq\mathbb W_3\subseteq\mathbb W_2\subseteq \mathbb W_1\subseteq \mathbb W_0\subseteq\mathbb L$ of the set of framed $m$-component links $\mathbb L=\mathbb L(m)$. Note that $\mathbb W_0$ consists of links that are evenly framed because a 
component has even framing if and only if it bounds a framed immersed 
disk in $B^4$.

In order to detect what stage of the filtration a particular link lies in, it would be convenient to define a set measuring the difference between $\mathbb W_{n}$ and $\mathbb W_{n+1}$. Because these are sets and not groups, the quotient is not defined. However we can still define an associated graded set in the following way.

Suppose $\mathcal W$ is an order $n+1$ Whitney tower in $M=S^3\times [0,1]$ where each of the order $0$ surfaces $A_1,\ldots,A_m$ is an annulus with one boundary component in $S^3\times\{0\}$ and one in $S^3\times\{1\}$. Then we say that the link $\partial_0\mathcal W$ is \emph{order $n+1$ Whitney tower concordant} to $\partial_1\mathcal W$. 
This allows us to define the associated graded set ${\sf W}_n$ as ${\mathbb W}_n$ modulo order $n+1$ Whitney tower concordance. Knots have a well-defined connected sum operation, but the analogous band-sum operation for links is not well-defined, even up to concordance. This makes the following proposition somewhat surprising; it follows from Theorem~\ref{thm:intro-order-raising}.

\begin {prop}[\cite{CST1}]
Band sum of links induces a well-defined operation which makes each ${\sf W}_n$ into a finitely generated abelian group.
\end {prop}

Our goal is to determine these groups ${\sf W}_n$.

\subsection{Free Lie and quasi-Lie algebras}\label{subsec:bracket-on-Lie}

Let $\sL=\sL(m)$ denote the free Lie algebra (over the ground ring $\Z$) on generators $\{X_1,X_2,\ldots,X_m\}$. 
It is $\N$-graded, $\sL=\oplus_n \sL_n$, where the degree~$n$ part $\sL_n$ is the
additive abelian group of length $n$ brackets, modulo Jacobi
relations and the self-annihilation relation $[X,X]=0$.
The free \emph{quasi-Lie algebra} ${\sf L}'$ is gotten from $\sL$ by replacing the self-annihilation relation with the
weaker anti-symmetry relation $[X,Y]=-[Y,X]$.

The bracketing map ${\sf L}_1\otimes {\sf L}_{n+1}\to {\sf L}_{n+2}$, has a nontrivial kernel, denoted ${\sf D}_n$. The analogous bracketing map on the free quasi-Lie algebra is denoted ${\sf D}'_n$. For later purposes, we now define a homomorphism $s\ell_{2n}\colon {\sf D}_{2n}\to \mathbb Z_2\otimes{\sf L}_{n+1}$. Given an element $X$, of ${\sf D}_{2n}$, its image under the bracketing map is zero in ${\sf L}_{2n+2}$. However, regarding the bracket as being in ${\sf L}'_{2n+2}$, we get an element of the kernel of the projection ${\sf L}'_{2n+2}\to{\sf L}_{2n+2}$. This kernel is isomorphic to $\mathbb Z_2\otimes{\sf L}_{n+1}$ by \cite{L3}, and so we get an element $s\ell_{2n}(x)$ of $\mathbb Z_2\otimes {\sf L}_{n+1}$ as desired.

\subsection{The total Milnor invariant}
Let $L$ be a link where all the
longitudes lie in $\Gamma_{n+1}$, the $(n+1)$th term of the
lower central series of the link group $\Gamma:=\pi_1(S^3\setminus
L)$.  By Van Kampen's Theorem $
\frac{\Gamma_{n+1}}{\Gamma_{n+2}}\cong\frac{F_{n+1}}{F_{n+2}}\cong {\sf L}_{n+1}
$,
where $F=F(m)$ is the free group on meridians. Let $ \mu_n^i(L)\in{\sf L}_{n+1}$ denote the image of the $i$-th longitude.
 The \emph{total Milnor invariant}
$\mu_n(L)$ of order~$n$ is defined by
\[
\mu_n(L):=\sum_i X_i \otimes \mu_n^i(L) \in {\sf L}_1\otimes {\sf L}_{n+1}
\]
It turns out that in fact $\mu_n(L)\in{\sf D}_n$ (by ``cyclic symmetry''). The invariant $\mu_n(L)$ is a convenient way of packaging all Milnor invariants of length $n+2$ in one piece.

\begin {thm} [\cite{CST4}]\label{thm:framed-mu-even-and-odd}
For all $n\in\mathbb N$, the total Milnor invariant is a well-defined homorphism $\mu_n\colon{{\mathsf W}}_n\to{\sf D}_n$ such that
\begin{enumerate}
\item  $\mu_n$ is an epimorphism for odd $n$; denote its kernel by ${\sf K}^\mu_{n}$.
\item For even $n$, $\mu_n$ is a monomorphism with image ${\sf D}'_n<{\sf D_n}$. 
\end{enumerate}
\end {thm}
So $\mu_{n}$ is an algebraic obstruction for $L$ bounding a Whitney tower of order $n+1$ which is a complete invariant in half the cases. In the other half, we'll need the following additional invariants.

\subsection{Higher-order Sato-Levine invariants}
Suppose $L\in\mathbb W_{2n-1}$ and $\mu_{2n-1}(L)=0$. This implies that the longitudes lie in $\Gamma_{2n}$, so that $\mu_{2n}(L)\in {\sf D}_{2n}$ is defined. Define the \emph{order $2n-1$ Sato-Levine invariant} by 
$\operatorname{SL}_{2n-1}(L)=s\ell_{2n}\circ\mu_{2n}(L)$,
where $s\ell_{2n}$ is defined above.

\begin {thm}[\cite{CST4}]\label{thm:SL}
 For all $n$, the Sato-Levine invariant gives a well-defined epimorphism $\operatorname{SL}_{2n-1}\colon {\sf K}^\mu_{2n-1}\twoheadrightarrow \mathbb Z_2\otimes {\sf L}_{n+1}$. Moreover, it is an isomorphism for even $n$. 
 \end {thm}
The case $\operatorname{SL}_1$ is the original Sato-Levine \cite{Sa} invariant of a 2-component classical link, and we describe in \cite{CST4} (and below) how the $\operatorname{SL}_{2n-1}$ are obstructions to ``untwisting'' an order $2n$ twisted Whitney tower.

\subsection{Higher-order Arf invariants}
We saw above that the structure of the groups ${\sf W}_n$ is completely determined for $n\equiv 0,2,3\!\mod 4$ by Milnor and higher-order Sato-Levine invariants. 

\begin {thm}[\cite{CST4}]\label{thm:framed-Arf}
Let ${\sf K}^{{\rm SL}}_{4k-3}$ be the kernel of $\operatorname{SL}_{4k-3}$. Then there is an epimorphism
$\alpha_{k}\colon\mathbb Z_2\otimes {\sf L}_{k}\twoheadrightarrow 
{\sf K}^{\operatorname{SL}}_{4k-3}$.
\end {thm}

\begin {conj}\label{arfconj}
$\alpha_{k}$ is an isomorphism.
\end {conj}
This conjecture is true when $k=1$, and indeed the inverse map $\alpha_1^{-1}\colon{\sf W}_1\to\mathbb Z_2\otimes{\sf L}_1$ is given by the classical Arf invariant of each component of the link. 

Regardless of whether or not Conjecture~\ref{arfconj} is true, $\alpha_{k}$ induces an isomorphism $\overline{\alpha}_{k}$ on $(\mathbb Z_2\otimes {\sf L}_{k})/\Ker \alpha_{k}$.

\begin  {defn}\label{def:higher-order-Arf}  
The \emph{higher-order Arf invariants} are defined by
$$
\operatorname{Arf}_{k}:=(\overline{\alpha}_{k})^{-1}:{\sf K}^{\operatorname{SL}}_{4k-3}\to(\mathbb Z_2\otimes {\sf L}_{k})/\Ker \alpha_{k}
$$
%to be the inverse of $\overline{\alpha}_{4n+1}$.
\end  {defn}
Any of the $\operatorname{Arf}_{k}$ which are non-trivial would the only possible remaining obstructions to a link bounding a  Whitney tower of order $4k-2$, following the Milnor and Sato-Levine invariants: 

\begin {cor}[\cite{CST4}]\label{cor:mu-sl-arf-classify} 
The associated graded groups ${\mathsf W}_n$ are classified by 
$\mu_n$, $\operatorname{SL}_n$ if $n$ is odd, and, for $n=4k-3$, $\operatorname{Arf}_k$.
\end {cor}

The first unknown Arf invariant is $\operatorname{Arf}_2\colon{\sf W}_5\to\mathbb  Z_2\otimes {\sf L}_2$, which in the case of $2$-component links would be a $\mathbb Z_2$-valued invariant, evaluating non-trivially on the Bing double of any knot with non-trivial classical Arf invariant. Evidence supporting the existence of non-trivial $\operatorname{Arf}_{k}$ is provided by the fact that such links are known to not be slice \cite{Cha}.
All cases for $k>1$ are currently unknown, but if $\operatorname{Arf}_2$ is trivial then all higher-order $\operatorname{Arf}_{k}$
would also be trivial \cite{CST2}.

%%%%%%%%%%%%%%%%%%%%%%%%%%%%%%%%%%%%%%%%%%%%

%%%%%%%%%%%%%%%%%%%%%%%%%%%%%%%%%%%%%%

\section{Twisted Whitney towers}
The order $n$ Sato-Levine invariants are defined as a certain projection of order $n+1$ Milnor invariants, suggesting that a slightly modified version of the Whitney tower filtration would put the Milnor invariants all in the right order, with no more need for the Sato-Levine invariants. In this section we discuss how this corresponds to the geometric notion of \emph{twisted} Whitney towers.

\begin{figure}[ht!]
        \centerline{\includegraphics[width=.5\columnwidth]{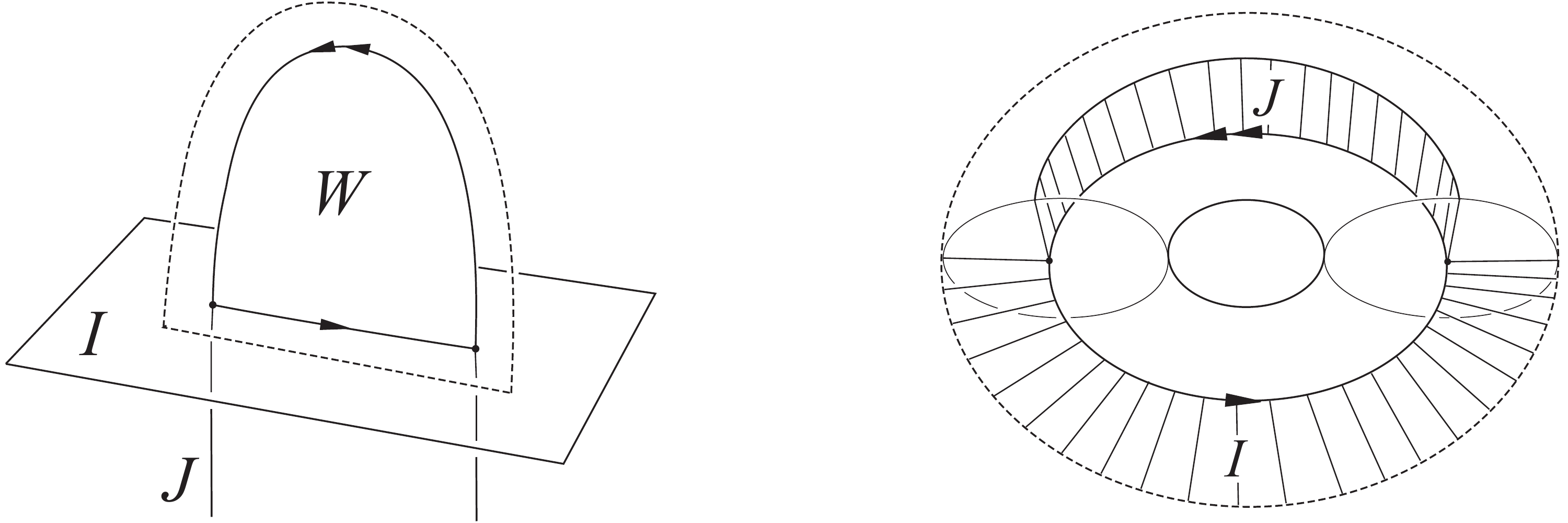}}
        \caption{The Whitney section over the boundary of a framed Whitney disk is
         indicated by the dotted loop. }
        \label{Framing-of-Wdisk-fig}
\end{figure}

%================
\subsection{Twisted Whitney disks}\label{subsec:twisted-w-disks}
The normal disk-bundle of a Whitney disk $W\looparrowright M$ is isomorphic to $D^2\times D^2$,
and comes equipped with a canonical nowhere-vanishing \emph{Whitney section} over the boundary given by pushing $\partial W$  tangentially along one sheet and normally along the other.

 In Figure~\ref{Framing-of-Wdisk-fig}, the Whitney section is indicated by a dotted loop shown on the left for a Whitney disk $W$ pairing intersections between surface sheets $I$ and $J$ in a $3$-dimensional slice of $4$--space. On the right is shown an embedding into $3$--space of the normal disk-bundle over $\partial W$, indicating how the Whitney section determines a nowhere-vanishing section which lies in the $I$-sheet and avoids the $J$-sheet.
 
%(see Figure~\ref{Framing-of-Wdisk-fig} and e.g.~\cite{FQ}).
The Whitney section determines the 
relative Euler number $\omega(W)\in\Z$ which represents the obstruction to extending
the Whitney section across $W$. It depends only on a choice of orientation of the tangent bundle of the ambient 4-manifold restricted to the Whitney disk, i.e.\ a local orientation. Following traditional terminology, when $\omega(W)$ vanishes $W$ is said to be \emph{framed}. (Since $D^2\times D^2$ has a unique trivialization up to homotopy, this terminology is only mildly abusive.)
If $\omega(W)=k$, we say that $W$ is
$k$-\emph{twisted}, or just \emph{twisted} if the value of $\omega(W)$ is not specified.

%%%%%%%%%%%%%%%%%%%%%%%%%%%%%%%%%%%%%%%%%%%
\begin{figure}[ht!]
        \centerline{\includegraphics[width=.7\columnwidth]{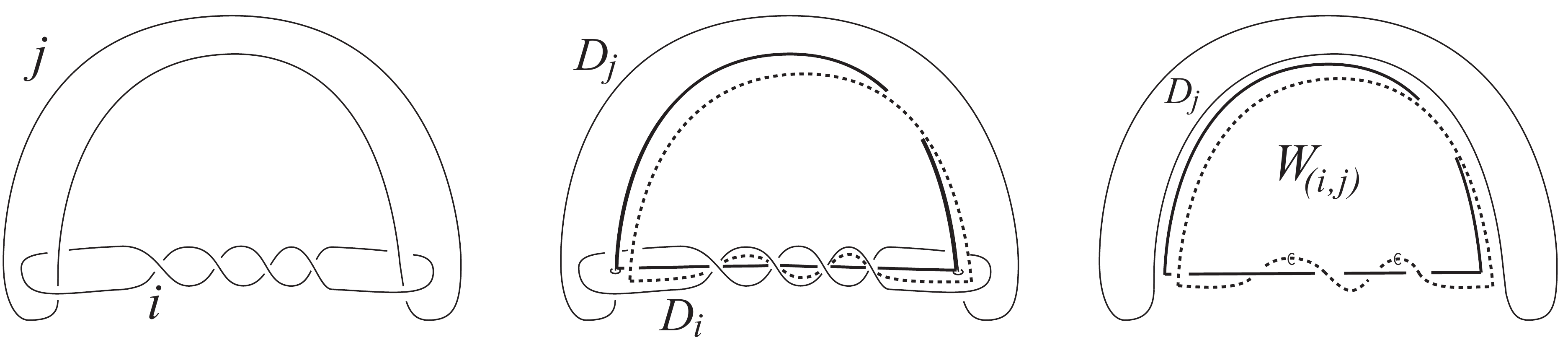}}
        \caption{Pushing into the $4$--ball from left to right: An $i$- and $j$-labeled twisted Bing double of the unknot bounds disks $D_i$ and $D_j$, which support a $2$-twisted Whitney disk $W_{(i,j)}$.
        The Whitney section is indicated by the dotted loop in the bottom center, and the intersections between its extension and the Whitney disk are shown in the bottom right.}
        \label{fig:Bing-unlink-W-disk}
\end{figure}
%%%%%%%%%%%%%%%%%%%%%%%%%%%%%%%%%%%%%%%%%%%%

In the definition of an order $n$ Whitney tower given above 
all Whitney disks are required to be framed ($0$-twisted). It turns out that the natural generalization to twisted Whitney towers involves allowing non-trivially twisted Whitney disks only in at least ``half the order'' as follows:

\begin  {defn}\label{def:twisted-W-towers}
A \emph{twisted Whitney tower of order $(2n-1)$} is just a (framed) Whitney
tower of order $(2n-1)$ as in Definition~\ref{def:w-tower} above.

A \emph{twisted Whitney tower of order $2n$} is a Whitney
tower having all intersections of order less than $2n$ paired by
Whitney disks, with all Whitney disks of order less than $n$ required to be framed, 
but Whitney disks of order at least $n$ allowed to be $k$-twisted for any $k$.
\end  {defn}

Note that, for any $n$, an order $n$ (framed) Whitney tower
is also an order $n$ twisted Whitney tower. We may
sometimes refer to a Whitney tower as a \emph{framed} Whitney tower to emphasize
the distinction, and will always use the adjective ``twisted'' in the setting of
Definition~\ref{def:twisted-W-towers}.

%\begin{rem}\label{rem:motivate-twisted-towers}
%The convention of allowing only order $\geq n$ twisted Whitney disks in order $2n$ twisted Whitney towers will be explained both algebraically and geometrically in \cite{CST2}. It turns out that any order $2n$ twisted Whitney tower can be modified
%so that all its Whitney disks of order $>n$ are framed, so the twisted Whitney disks of order equal to $n$ are the important ones.  
%\end{rem}
%=================
\subsection{Twisted Whitney tower concordance}\label{subsec:twisted-w-conc}

Let $\mathbb W^\iinfty_{n}$ be the set of framed links in $S^3$ which are boundaries of order $n$ twisted Whitney towers in $B^4$,
with no requirement that the link framing is induced by the order $0$ disks.  Notice that $\mathbb W^\iinfty_{2n-1}=\mathbb W_{2n-1}$. 
 While not immediately obvious, it is true that this defines a filtration $\cdots\subseteq\mathbb W^\iinfty_3\subseteq\mathbb W^\iinfty_2\subseteq \mathbb W^\iinfty_1\subseteq \mathbb W^\iinfty_0=\mathbb L$. 
 As in the framed setting above, letting ${\sf W}^\iinfty_n$ be the set $\mathbb W^\iinfty_n$ modulo order $(n+1)$ twisted Whitney tower concordance yields a finitely generated abelian group.

\begin {thm}[\cite{CST2,CST4}]\label{thm:mu-epimorphism-on-W-infty}
The total Milnor invariants give epimorphisms $\mu_n\colon{\sf W}_n^\iinfty\twoheadrightarrow {\sf D}_n$ which are isomorphisms for $n\equiv 0,1,3\!\mod 4$. Moreover, the kernel ${\sf K}^\iinfty_{4k-2}$ of $\mu_{4k-2}$ is isomorphic to the kernel ${\sf K}^{\rm SL}_{4k-3}$ of the Sato-Levine map from the previous section.
 \end {thm}
 
Conjecture~\ref{arfconj} hence says that ${\sf K}^\iinfty_{4k-2} \cong \Z_2\otimes {\sf L}_{k}$ and our Arf-invariants $\operatorname{Arf}_{k}$ represent the only remaining obstruction to a link bounding an order $4k-1$ twisted Whitney tower:

 \begin {cor}\label{cor:twisted-classification}
 The groups ${\sf W}^\iinfty_{n}$ are classified by $\mu_n$ and, for $n=4k-2$,  $\operatorname{Arf}_k$.
 \end {cor}  
  
\subsection{Gropes and k-slice links}
Roughly speaking, a link is said to be ``$k$-slice'' if it is the boundary of a surface
which ``looks like a collection of slice disks modulo $k$-fold commutators in the fundamental group of the complement of the surface''.
Precisely, $L\subset S^3$ is \emph{$k$-slice} if $L$ bounds an embedded orientable surface $\Sigma\subset B^4$ such that $\pi_0(L)\to\pi_0(\Sigma)$ is a bijection and there is a push-off homomorphism $\pi_1(\Sigma)\to\pi_1(B^4\setminus\Sigma)$ whose image lies in the $k$th term of the lower central series 
$(\pi_1(B^4\setminus\Sigma))_k$. Igusa and Orr proved the following ``$k$-slice conjecture'' in \cite{IO}:

\begin {thm}[\cite{IO}]
A link $L$ is $k$-slice if and only if $\mu_{i}(L)=0$ for all $i\leq 2k-2$.
\end {thm}

A $k$-fold commutator in $\pi_1X$ has a nice topological model in terms of a continuous map $G\to X$, where $G$ is a \emph{grope of class $k$}. Such 2-complexes $G$ (with specified ``boundary'' circle) are recursively defined as follows.
 A grope of class $1$ is a circle. A grope of class $2$ is an orientable surface with one boundary component. A grope of class $k$ is formed by attaching to every dual pair of basis curves on a class $2$ grope
 a pair of gropes whose classes add to $k$. A curve $\gamma\colon S^1\to X$ in a topological space $X$ is a $k$-fold commutator if and only if it extends to a continuous map of a grope of class $k$. Thus one can ask whether being $k$-slice implies there is a basis of curves on $\Sigma$ that bound \emph{disjointly embedded} gropes of class $k$ in $B^4\setminus \Sigma$. Call such a link \emph{geometrically $k$-slice}. 
 
 \begin {prop}[\cite{CST2}]\label{prop:geo-k-slide-equals-odd-w-tower}
 A link $L$ is geometrically $k$-slice if and only if $L\in\mathbb W^\iinfty_{2k-1}$. 
 \end {prop}
 
 This is proven using a construction from \cite{S} which allows one to freely pass between class $n$ gropes and order $n-1$ Whitney towers.
 So the higher-order Arf invariants $\operatorname{Arf}_k$ detect the difference between $k$-sliceness and geometric $k$-sliceness. It turns out that every $\operatorname{Arf}_k$ value can be realized by (internal) band summing iterated Bing doubles of the figure-eight knot. Every Bing double is a boundary link, and one can choose the bands so that the sum remains a boundary link. 
 This implies

\begin {thm}[\cite{CST2}]
 A link $L$ has vanishing Milnor invariants of orders up to $2k-2$ if and only if it is geometrically $k$-slice after connected sums with internal band sums of iterated Bing doubles of the figure-eight knot.
 \end {thm}
 
Here (and in Theorem~\ref{thm:Artin-kernel} below), the figure-eight knot can be replaced by any knot with non-trivial (classical) Arf invariant.

The added boundary links in the above theorem bound disjoint surfaces in $S^3$ which clearly allow immersed disks in $B^4$ bounded by curves representing a basis of first homology. 
In \cite{CST2} we will show that this implies:

\begin {cor}
A link  has vanishing Milnor invariants of orders up to $2k-2$ if and only if its components bound disjointly embedded surfaces $\Sigma_i\subset B^4$, with each surface a connected sum of two surfaces $\Sigma'_i$ and $\Sigma''_i$ such that
\begin{enumerate}
\item
 a basis of curves on $\Sigma'_i$ bound disjointly embedded framed gropes of class $k$ in the complement of $\Sigma:= \cup_i\Sigma_i$,
 \item a basis of curves on $\Sigma''_i$ bound immersed disks in the
 complement of $\Sigma\cup G$, where $G$ is the union of all class $k$ gropes on the surfaces $\Sigma_i'$.
\end{enumerate}
\end {cor}

This is an enormous geometric strengthening of Igusa and Orr's result who, under the same assumption on the vanishing of Milnor invariants, show the existence of a surface $\Sigma$ and maps of class $k$ gropes, with no control on their intersections and self-intersections. Our proof uses the full power of the obstruction theory for Whitney towers, whereas they do a sophisticated computation of the third homology of the groups $F/F_{2k}$.

 \subsection{String links and the Artin representation}
Let $L$ be a string link with $m$ strands embedded in $D^2\times[0,1]$. By Stallings' Theorem \cite{Sta}, the inclusions 
$(D^2\setminus \{\mbox{$m$ points}\})\times\{i\}\hookrightarrow (D^2\times[0,1])\setminus L$ for $i=0,1$ induce isomorphisms on all lower central quotients of the fundamental groups. In fact, the induced automorphism of the lower central quotients 
$F/F_n$ of the free group $F=\pi_1(D^2\setminus \{\mbox{$m$ points}\})$ is explicitly characterized by conjugating the meridional generators of $F$ by longitudes. Let $\operatorname{Aut}_0(F/F_n)$ consist of those automorphisms of $F/F_n$ which are defined by conjugating each generator and which fix the product of generators. This leads to the \emph{ Artin representation} $\mathbb{SL}\to \operatorname{Aut}_0(F/F_{n+2})$ where $\mathbb{SL}$ is the set of concordance classes of pure framed string links.

The set of string links has an advantage over links in that it has a well-defined monoid structure given by stacking. Indeed, modulo concordance, it becomes a (noncommutative) group. Whitney tower filtrations can also be defined in this context, giving rise to filtrations $\mathbb{SW}_n$ and $\mathbb{SW}_n^\iinfty$ of this group $\mathbb{SL}$.

\begin {thm}[\cite{CST5}] \label{thm:Artin-kernel}
The sets $\mathbb{SW}_n$ and $\mathbb{SW}_n^\iinfty$ are normal subgroups of $\mathbb{SL}$ which are central modulo the next order. We obtain nilpotent groups $\mathbb {SL}/\mathbb{SW}_n$ and
$\mathbb {SL}/\mathbb{SW}^\iinfty_n$ and the associated graded are isomorphic to our previously defined groups
\[
\W_n \cong\mathbb {SW}_n/\mathbb{SW}_{n+1}  \quad \text{ and } \quad \mathbb {SW}^\iinfty_n/\mathbb{SW}^\iinfty_{n+1} \cong \W^\iinfty_n
\]
Finally, the Artin representation induces a well-defined epimorphism 
$\operatorname{Artin}_{n}\colon\mathbb {SL}/\mathbb{SW}^\iinfty_{n}\twoheadrightarrow\operatorname{Aut}_0(F/F_{n+2})$ whose kernel is generated by internal band sums of iterated Bing doubles of the figure-eight knot.
\end {thm}

The Artin representation is thus an invariant on the whole group $\mathbb {SL}/\mathbb{SW}^\iinfty_n$, not just on the associated graded groups as in the case of links. It packages the total Milnor invariants $\mu_k, k=0,
\dots, n$ on string links together into a group homomorphism.

%\begin {thm}[\cite{CST5}]
%The following subsets of $\mathbb {SL}/\mathbb{SW}^\iinfty_n$ are equal: 
%\begin{enumerate}
%\item The subgroup generated by iterated Bing doubles of knots with non-trivial (classical) Arf invariant.
%\item The set of boundary string links.
%\item The set of $\pi_1$-null string links.
%\end{enumerate}
%Moreover, this subgroup (determined by any of these three definitions) is exactly the kernel of the Artin map $\mathbb {SL}/\mathbb{SW}^\iinfty_n\to\operatorname{Aut}_0(F/F_{n+2})$.
%\end {thm}
%Here a \emph{boundary string link} is a string link whose standard closure along $\partial(D^2\times[0,1])$ bounds disjoint surfaces in $D^2\times[0,1]$; and a 
%\emph{$\pi_1$-null string link} $L$ is a string link whose standard closure bounds a surface $\Sigma$ in $M=D^2\times[0,1]\times[0,1]$ such that $\pi_0(L)\to\pi_0(\Sigma)$ is a bijection and
%which for which there is a push-off inducing the trivial homomorphism $\pi_1(\Sigma)\to\pi_1(B^4\setminus\Sigma)$.

%=================================================================================================================================

\section{Higher-order intersection invariants}\label{sec:geo-ints}
Proofs of the above results depend on two essential ideas: The higher-order intersection theory
of Whitney towers comes with an obstruction theory whose associated invariants take values in abelian groups of (un-rooted) trivalent trees. And by mapping to rooted trees, which correspond to iterated commutators, the obstruction theory for Whitney towers in $B^4$ can be identified with
algebraic invariants of the bounding link in $S^3$. A critical connection between these ideas is provided by the resolution of the Levine Conjecture (see below), which says that this map is an isomorphism.

\begin{figure}[ht!] 
        \centerline{\includegraphics[width=.7\columnwidth]{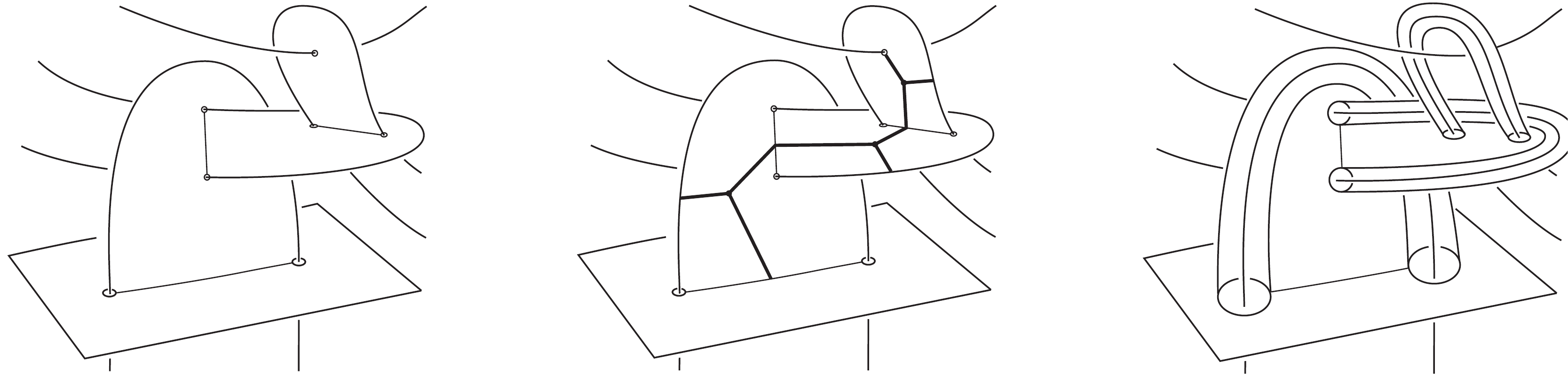}}
        \caption{From left to right: An unpaired intersection in a Whitney tower, (part of) its associated tree, and the result of surgering to a grope.}
        \label{fig:order3Whitneytower-and-withTrees-and-grope}
\end{figure}

In fact, it can be arranged that all singularities in a Whitney tower are contained in $4$-ball neighborhoods of the associated trivalent trees, which sit as embedded `spines';
and all relations among trees in the target group are realized by controlled manipulations of the Whitney disks. Mapping to rooted trees
corresponds geometrically to surgering Whitney towers to gropes, and these determine iterated commutators
of meridians of the Whitney tower boundaries as in Figure~\ref{fig:order3Whitneytower-and-withTrees-and-grope}.

\subsection{Trees and intersections}\label{subsec:trees-and-ints}
All trees are unitrivalent, with cyclic orderings of the edges at all trivalent vertices, and univalent vertices labeled from an index set $\{1,2,3,\ldots,m\}$. A \emph{rooted} tree has one unlabeled univalent vertex designated as the \emph{root}. Such rooted trees correspond to formal non-associative bracketings of elements from
the index set. The \emph{rooted product} $(I,J)$ of rooted trees $I$ and $J$ is the rooted tree gotten
by identifying the root vertices of $I$ and $J$ to a single vertex $v$ and sprouting a new rooted edge at $v$.
This operation corresponds to the formal bracket, and we identify rooted trees with formal brackets.
The \emph{inner product}  $\langle I,J \rangle $ of rooted trees $I$ and $J$ is the
unrooted tree gotten by identifying the roots of $I$ and $J$ to a single non-vertex point.
Note that 
all the univalent vertices of $\langle I,J \rangle $ are labeled.

The \emph{order} of a tree, rooted or unrooted, is defined to be the number of trivalent vertices,
and the following associations of trees to Whitney disks and intersection points respects the notion
of order given in Definition~\ref{def:w-tower}.

%\subsection{Whitney towers and higher-order intersections}%\label{subsec:order-zero-w-towers-and-ints}
%A collection $A_1,\ldots,A_m\looparrowright (M,\partial M)$ of 
%connected surfaces in a $4$--manifold $M$ is a \emph{Whitney tower of order zero} if the $A_i$ are \emph{properly immersed} in the sense that the boundary is embedded in $\partial M$ and the interior is generically immersed in $M \smallsetminus \partial M$. 

%If $\partial A_i\hookrightarrow \partial M$ comes equipped with a non-vanishing normal vector field, then we require it to extend to a nonvanishing normal vector field on $A_i$. And if $\partial A_i$ is equipped with an orientation, then we require that this orientation is induced by an orientation of $A_i$. 

%A {\em framing} of $\partial A_i$ (respectively $A_i$) is by definition a trivialization of the normal bundle of the immersion. If $M$ is oriented, this is equivalent to an orientation and a non-vanishing normal vector field on $\partial A_i$ (respectively $A_i$).

%For $M=B^4$, this means that the self-intersections of $A_i$ come in canceling pairs.
 
%like the {\em zero-framing} of a knot in $S^3$...

To each order zero surface $A_i$ is associated
the order zero rooted tree consisting of an edge with one vertex labeled by $i$, and
to each transverse intersection $p\in A_i\cap A_j$ is associated the order zero
tree $t_p:=\langle i,j \rangle$ consisting of an edge with vertices labeled by $i$ and $j$. 
The order 1 rooted Y-tree $(i,j)$, with a single trivalent vertex and two univalent labels $i$ and $j$,
is associated to any Whitney disk $W_{(i,j)}$ pairing intersections between $A_i$ and $A_j$. This rooted tree
can be thought of as an embedded subset of $M$, with its trivalent vertex and rooted
edge sitting in $W_{(i,j)}$, and its two other edges descending into $A_i$ and $A_j$ as sheet-changing paths. 

\begin{figure}[ht!]
        \centerline{\includegraphics[width=.5\columnwidth]{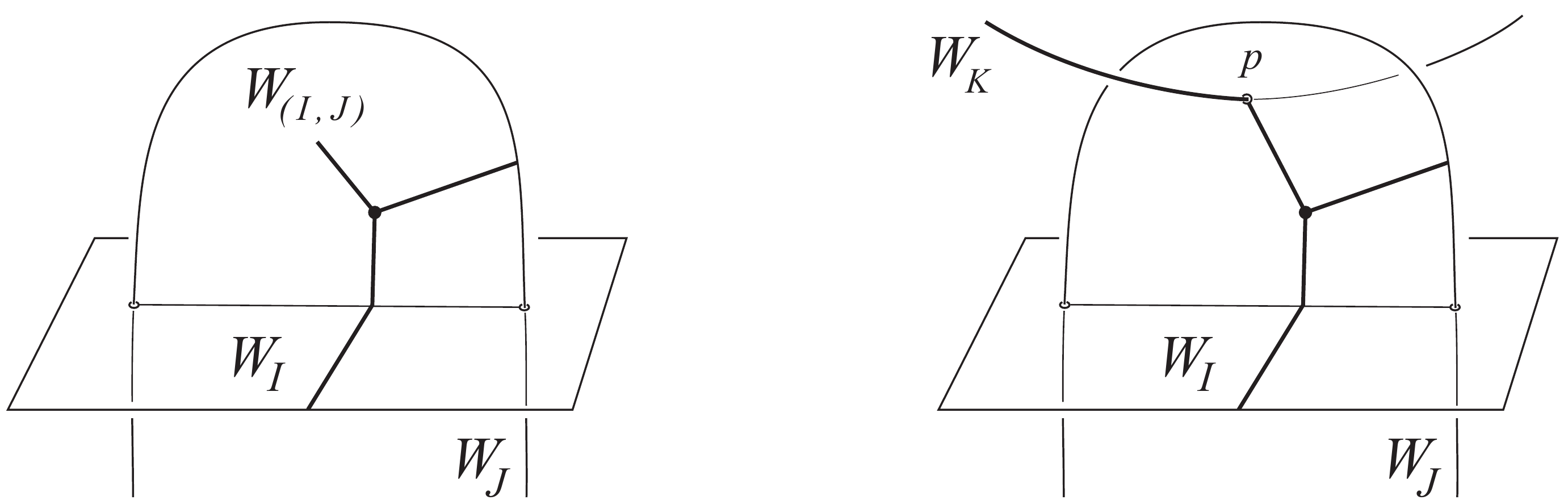}}
        \caption{
%On the left, (part of) the rooted tree $(I,J)$ associated to a Whitney disk $W_{(I,J)}$. On the right, (part of) the unrooted tree $t_p=\langle (I,J),K \rangle$ associated to an intersection $p\in W_{(I,J)}\cap W_K$. Note that $p$ corresponds to where the roots of $(I,J)$ and $K$ are identified to a (non-vertex) point in $\langle (I,J),K \rangle$.
}
        \label{WdiskIJandIJKint-fig}
\end{figure}

Recursively, the rooted tree $(I,J)$ is associated to any Whitney disk $W_{(I,J)}$ pairing intersections
between $W_I$ and $W_J$
(see left-hand side of Figure~\ref{WdiskIJandIJKint-fig}); 
with the understanding that if, say, $I$ is just a singleton $i$, then $W_I$ denotes the order zero surface $A_i$. 
To any transverse intersection $p\in W_{(I,J)}\cap W_K$ between $W_{(I,J)}$ and any
$W_K$ is associated the un-rooted tree $t_p:=\langle (I,J),K \rangle$  
(see right-hand side of Figure~\ref{WdiskIJandIJKint-fig}).

%%====================================
%  

\subsection{Intersections trees for Whitney towers}\label{subsec:intro-w-tower-int-trees}
The group $\cT_n$ (for each $n=0,1,2\ldots$) is the free abelian group on (unitrivalent labeled vertex-oriented)
order $n$ trees, modulo the usual AS (antisymmetry) and IHX (Jacobi) relations:
\begin{figure}[h]
\centerline{\includegraphics[width=.5\columnwidth]{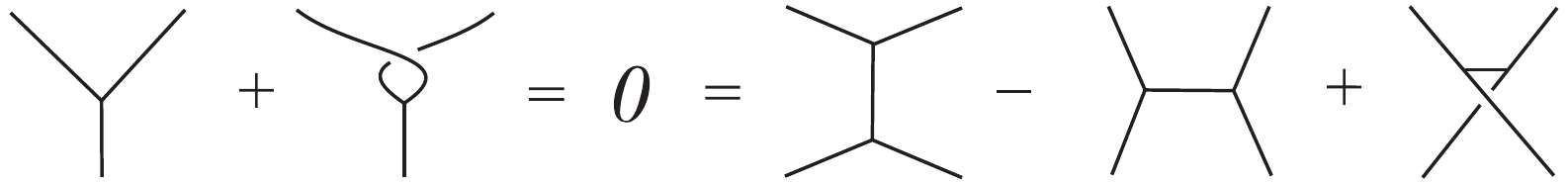}}
%        \caption{Local pictures of the \emph{antisymmetry} (AS) and \emph{Jacobi} (IHX) relations
%         in $\cT_n$. Here all trivalent orientations are induced from a fixed orientation of the plane, and univalent vertices possibly extend to subtrees which are fixed in each equation.}
         \label{fig:ASandIHXtree-relations}
\end{figure}

In even orders we define $\widetilde{\cT}_{2n}:=\cT_{2n}$, and in odd orders $\widetilde{\cT}_{2n-1}$ is defined to be the quotient of $\cT_{2n-1}$ by the \emph{framing relations}. 
These framing relations are defined as the image of homomorphisms $\Delta_{2n-1}:\Z_2\otimes\cT_{n-1}\rightarrow \cT_{2n-1}$ which are defined for generators $t\in\cT_{n-1}$ by
$\Delta (t):=\sum_{v\in t} \langle i(v),(T_v(t),T_v(t))\rangle$,
where $T_v(t)$ denotes the rooted tree gotten by replacing $v$ with a root, and the sum is over all univalent vertices of $t$, with $i(v)$ the original label of the univalent vertex $v$. 
%That $\Delta_{2n-1}$ is well-defined as a homomorphism on $\cT_{n-1}$ is clear since AS and IHX relations go to doubled relations. 
%The image of $\Delta_{2n-1}$ is $2$-torsion by AS relations and hence it factors through $\Z_2\otimes\cT_{n-1}$.
%See Figure~\ref{fig:framing-relations-order-one-and-three} for explicit illustrations of the framing relations in orders $1$ and $3$.
For example, in orders 1 and 3 the framing relations are:
\begin{figure}[h]
\centerline{\includegraphics[width=.5\columnwidth]{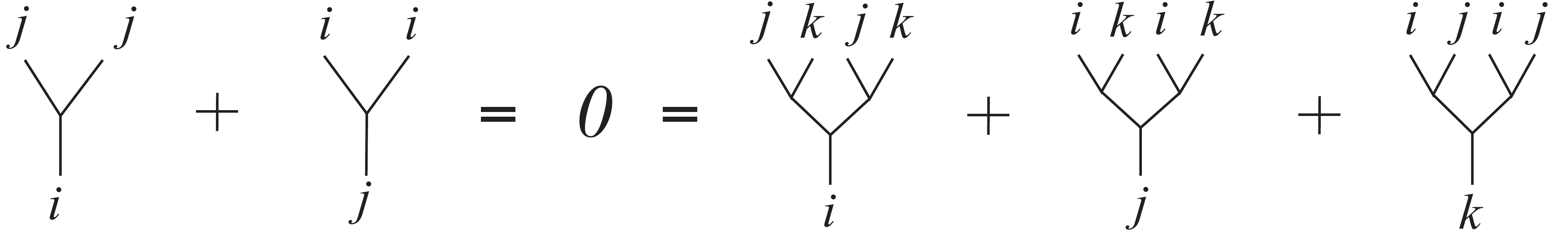}}
         %\caption{The order $1$ and order $3$ framing relations.}
         \label{fig:framing-relations-order-one-and-three}
\end{figure}

The obstruction theory works as follows:

\begin  {defn}
The \emph{order $n$ intersection tree} $\tau_n(\cW)$ of an order
$n$ Whitney tower $\cW$ is defined to be
$$
\tau_n(\cW):=\sum \epsilon_p\cdot t_p \in\widetilde{\cT}_n
$$ 
where the sum is over all order $n$ intersections $p$,
with
$\epsilon_p=\pm 1$ the usual sign of a transverse intersection
point (via certain orientation conventions, see e.g.~\cite{CST1}).
\end  {defn}

All relations in $\widetilde{\cT}_n$ can be realized by controlled
manipulations of Whitney towers, and further maneuvers allow 
algebraically canceling pairs of tree generators to be converted into intersection-point pairs admitting Whitney disks.
As a result, we get the following partial recovery of the ``algebraic cancelation implies geometric cancellation''
principle available in higher dimensions:  

\begin {thm}[\cite{CST1}]\label{thm:framed-order-raising-on-A}
If a collection $A$ of properly immersed surfaces in a simply-connected $4$--manifold supports an order $n$ Whitney tower 
$\cW$ with $\tau_n(\cW)=0\in\widetilde{\cT}_n$, then $A$ is homotopic (rel $\partial$) to 
$A'$ which supports an order $n+1$ Whitney tower.
\end {thm}

%

%

%

%%========================================

\subsection{Intersections trees for twisted Whitney towers}\label{subsec:intro-twisted-w-tower-int-trees}

For any rooted tree $J$ we define the corresponding {\em $\iinfty$-tree} (``twisted-tree''), denoted by $J^\iinfty$, by labeling the root univalent vertex with the symbol ``$\iinfty$'' (which will represent a ``twist'' in a Whitney disk normal bundle):
$J^\iinfty := \iinfty\!-\!\!\!- J $.

%The intersection theory for Whitney towers is extended to
%twisted Whitney towers as follows (details in \cite{CST1}):

%%
\begin  {defn}
The abelian group $\cT^{\iinfty}_{2n-1}$ is the quotient of $\widetilde{\cT}_{2n-1}$ by the \emph{boundary-twist relations}: 
\[
\langle  (i,J),J \rangle \,=\, i\,-\!\!\!\!\!-\!\!\!<^{\,J}_{\,J}\,\,=\,0
\] 
Here $J$ ranges over all order $n-1$ rooted trees (and the first equality is just a reminder of notation).

%This is the same as taking the quotient of $\widetilde{\cT}_{2n-1}=\cT_{2n-1}/\im(\Delta_{2n-1})$ by boundary-twist relations since $\im(\Delta_{2n-1})$ is contained in the span of boundary-twist relations -- see Section~\ref{sec:proof-thm-odd}.

%
%The boundary-twist relations
%correspond geometrically to the fact that 
%performing a boundary twist (Figure~\ref{boundary-twist-and-section-fig}) on an order $n$ Whitney disk $W_{(i,J)}$ creates an order $2n-1$ intersection point
%$p\in W_{(i,J)}\cap W_J$ with associated tree $t_p=\langle  (i,J),J \rangle $ (which is 2-torsion
%by the AS relations) and changes $\omega (W_{(i,J))})$ by $\pm1$. Since order $n$ twisted Whitney disks are allowed in an order $2n$ Whitney tower such trees do not represent obstructions to the existence of the next order twisted tower.

The abelian group $\cT^{\iinfty}_{2n}$ is gotten from $\widetilde{\cT}_{2n}=\cT_{2n}$ by including order $n$ $\iinfty$-trees as new generators
and introducing the 
following new relations (in addition to the IHX and antisymmetry relations on non-$\iinfty$ trees):
$$
J^\iinfty=(-J)^\iinfty\quad\quad
 I^\iinfty=H^\iinfty+X^\iinfty- \langle H,X\rangle \quad\quad
 2\cdot J^\iinfty=\langle J,J\rangle 
 $$
\end  {defn}
 The left-hand \emph{symmetry} relation corresponds to the fact that the framing obstruction on a Whitney disk is independent of its orientation; the middle \emph{twisted IHX} relations can be realized by a Whitney move near a twisted Whitney disk, and the right-hand \emph{interior twist} relations can be realized by cusp-homotopies in Whitney disk interiors. As described in \cite{CST4}, the twisted groups $\cT^{\iinfty}_{2n}$ can naturally be identified with a universal quadratic refinement of the $\cT_{2n}$-valued intersection pairing $\langle \, \cdot,\cdot \,\rangle$ on framed Whitney disks. 

Recalling from Definition~\ref{def:twisted-W-towers} that twisted Whitney disks
only occur in even order twisted Whitney towers, intersection trees for twisted Whitney towers
are defined as follows:

\begin  {defn}\label{def:tau-infty}
The \emph{order $n$ intersection tree}
$\tau_{n}^{\iinfty}(\cW)$ of an order
$n$ twisted Whitney tower $\cW$ is defined to be
$$
\tau_{n}^{\iinfty}(\cW):=\sum \epsilon_p\cdot t_p + \sum \omega(W_J)\cdot J^\iinfty\in\cT^{\iinfty}_{n}
$$
where the first sum is over all order $n$ intersections $p$ and the second sum is over all order $n/2$
Whitney disks $W_J$ with twisting $\omega(W_J)\in\Z$ (computed from a consistent choice of local orientations). 
\end  {defn}

 By ``splitting'' the twisted Whitney disks \cite{CST1}
it can be arranged that $|\omega(W_J)|\leq 1$, leading to signs like $\epsilon_p$ (or zero coefficients).
The obstruction theory also holds for twisted Whitney towers:

\begin {thm}[\cite{CST1}]\label{thm:twisted-order-raising-on-A}
If a collection $A$ of properly immersed surfaces in a simply-connected $4$--manifold supports an order $n$ twisted Whitney tower $\cW$ with $\tau_n^\iinfty(\cW)=0\in\cT^\iinfty_n$, then $A$ is homotopic (rel $\partial$) to 
$A'$ which supports an order $n+1$ twisted Whitney tower.
\end {thm}

%\begin{figure}[h]
%\centerline{\includegraphics[width=80mm]{twisted-Wdisk-on-figure-eight-knot}}
%         \caption{Realizing the $\iinfty$-tree $(i,j)^\iinfty$.}
%         \label{fig:twisted-Wdisk-on-figure-eight-knot}
%\end{figure}

\subsection{Remark on the framing relations}\label{rem:framing-relations} The framing relations 
in the untwisted groups $\widetilde{\cT}_{2n-1}$ correspond to the twisted IHX relations among $\iinfty$-trees in  
$\cT^{\iinfty}_{2n}$ via a geometric boundary-twist operation which converts an order $n$ $\iinfty$-tree $(i,J)^\iinfty$ to an order 
$2n-1$ (untwisted) tree $\langle  (i,J),J \rangle$.

\subsection{Realization maps}

\begin{figure}[h]
\centerline{\includegraphics[width=.5\columnwidth]{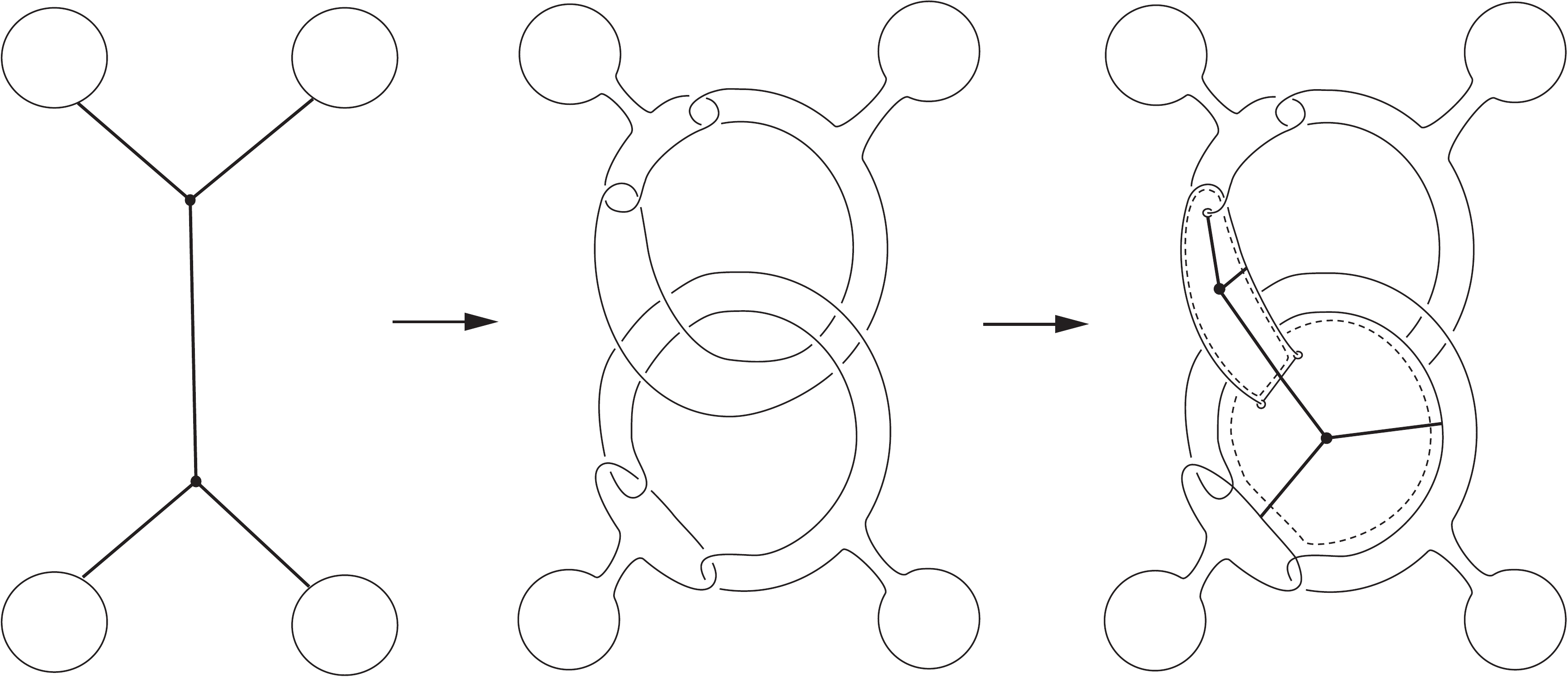}}
       \caption{Realizing an order $2$ tree in a Whitney tower by Bing doubling.}
         \label{fig:Cochran-Habiro-fig}
\end{figure}

In \cite{CST1} we describe how to construct surjective realization maps
$\widetilde R_n : \widetilde\cT_n \sra\W_n $ and
$ R^\iinfty_n : \cT^\iinfty_n \sra\W^\iinfty_n $ by applying the operation of iterated Bing doubling.
This construction is essentially the same as the well-known application of Habiro's clasper-surgery \cite{H} and Cochran's
realization of Milnor invariants \cite{C}, extended to twisted Bing doubling (Figures~\ref{fig:Cochran-Habiro-fig} and~\ref{fig:Bing-unlink-W-disk}). To prove the realization maps are well-defined we need to use Theorems
~\ref{thm:framed-order-raising-on-A} and ~\ref{thm:twisted-order-raising-on-A} respectively.

The above Conjecture~\ref{arfconj} on the non-triviality of the higher-order Arf invariants can be
succinctly rephrased as the assertion that the realization maps $\widetilde R_n$ and
$ R^\iinfty_n$ are \emph{isomorphisms} for all $n$.
Progress towards confirming this assertion -- namely complete answers in $3/4$ of the cases and partial answers in the remaining cases, as described by the above-stated results -- has been accomplished by identifying
intersection trees with Milnor invariants, as we describe next.

\subsection{Intersection trees and Milnor's link invariants}
%The set of $m$-component framed links in $S^3$ which bound order $n$ twisted Whitney towers in $B^4$ is 
%denoted $\bW^\iinfty_n$; and the quotient of $\bW^\iinfty_n$ by the equivalence relation of order $n+1$ twisted Whitney 
%tower concordance (immersed annuli in $S^3\times I$ supporting order $n+1$ twisted Whitney 
%towers) is denoted by $\W^\iinfty_n$.  For brevity we concentrate here on this twisted Whitney tower filtration, but all results
%carry over to framed Whitney tower and grope filtrations of link concordance (with some index-shifting -- see \cite{CST1}).

%As a result of the obstruction theory, $\W^\iinfty_n$ forms a group under connected sum; and a realization map 
%$R^\iinfty_n : \cT^\iinfty_n \to\W^\iinfty_n$ can be constructed using iterated Bing doublings of Hopf links and knots with non-trivial Arf invariant \cite{CST1}.

%Denoting by $\sL_k$ the degree part of the free $\Z$-Lie algebra on generators $\{X_1,X_2,\ldots,X_m\}$,
%the first non-vanishing order $n$ total Milnor invariant $\mu_n(L)$ of a link $L\subset S^3$ is defined by 
%$\mu_n(L):=\sum_i\,X_i\otimes\mu^i_n(L)\in\sL_1\otimes\sL_{n+1}$, where $\mu^i_n(L)\in\sL_{n+1}$ is determined by the $i$th longitude of $L$ (as an iterated commutator of meridians) via the identification of $\sL$ with the lower central series of the free group. It turns out that $\mu_n(L)$ lies in $\sD_n$, the kernel of the bracketing map $\sL_1\otimes\sL_{n+1}\to\sL_{n+2}$ \cite{HM,O}. (The Lie algebra $\sD_n$ also appears in a variety of other topological settings, e.g.~\cite{CV1,CV2,J,Kon,L1,Mor2}.)

The connection between intersection trees and Milnor invariants is via a surjective map 
$\eta_n:\cT^\iinfty_n \to \sD_n$ which converts trees to rooted trees (interpreted as Lie brackets) by summing over all ways of choosing a root:

For $v$ a univalent vertex of an order $n$ (un-rooted non-$\iinfty$) tree
%$t$ with univalent vertices labeled from the index set $\{1,2,\ldots,m\}$ and cyclic edge orderings at each trivalent vertex, 
denote by $B_v(t)\in\sL_{n+1}$ the
Lie bracket of generators $X_1,X_2,\ldots,X_m$ determined by the formal bracketing of indices
which is gotten by considering $v$ to be a root of $t$.

Denoting the label of a univalent vertex $v$ by $\ell(v)\in\{1,2,\ldots,m\}$, the
map $\eta_n:\cT^\iinfty_n\rightarrow \sL_1 \otimes \sL_{n+1}$
is defined on generators by
$$
\eta_n(t):=\sum_{v\in t} X_{\ell(v)}\otimes B_v(t)
\quad \, \,
\mbox{and}
\quad \, \,
\eta_n(J^\iinfty):=\frac{1}{2}\eta_n(\langle J,J \rangle)
$$
where the first sum is over all univalent vertices $v$ of $t$, and the second expression lies in $\sL_1 \otimes \sL_{n+1}$ 
because the coefficient of $\eta_n(\langle J,J \rangle)$ is even.

The proof of the following theorem (which implies Theorem~\ref{thm:mu-epimorphism-on-W-infty} above) shows that the map $\eta$ corresponds to a construction which converts Whitney towers into embedded gropes \cite{S}, via the grope duality of \cite{KT}:

\begin {thm}[\cite{CST2}]\label{thm:mu-equals-eta-of-tau}
If $L$ bounds a twisted Whitney tower $\cW$ of order $n$, then the total Milnor invariants $\mu_k(L)$ vanish for $k<n$, and
$\mu_n(L) = \eta_n \circ\tau^\iinfty_n(\cW) \in \sD_n$.
\end {thm}

Thus one needs to understand the kernel of $\eta_n$ before the obstruction theory can proceed. This is accomplished
by resolving \cite{CST3} a closely related conjecture of J. Levine \cite{L2}, as discussed next. 

\subsection{The Levine Conjecture and its implications}\label{subsec:Levine-conj}
The bracket map kernel $\sD_n$ turns out to be relevant to a variety of topological settings (see e.g.~the introduction to \cite{CST3}), and was known to be isomorphic to 
$\cT_n$ after tensoring with $\Q$, when Levine's study of the cobordism groups of $3$-dimensional homology cylinders \cite{L1,L2} led him to conjecture that $\cT_n$ is in fact isomorphic
to the quasi-Lie bracket map kernel $\sD_n'$, via the analogous map $\eta'_n$
which sums over all choices of roots (as in the left formula for $\eta$ above).

Levine made progress in \cite{L2,L3}, and in \cite{CST3} we affirm his conjecture:

\begin {thm}[\cite{CST3}]\label{thm:LC}

$\eta'_n:\cT_n\to\sD'_n$ is an isomorphism for all $n$.

\end {thm}
The proof of Theorem~\ref{thm:LC} uses techniques from discrete Morse theory on chain complexes, including an extension of the theory to complexes containing torsion.
A key idea involves defining combinatorial vector fields that are inspired by the Hall basis algorithm for free Lie algebras and its generalization by Levine to quasi-Lie algebras.

As described in \cite{CST4}, Theorem~\ref{thm:LC} has several direct applications to Whitney towers, including the completion of
the calculation of $\W^\iinfty_n$ in three out of four cases:

\begin {thm}[\cite{CST4}] \label{thm:Milnor-twisted} $\eta_n:\cT^\iinfty_n \to \sD_n$ are isomorphisms for $n\equiv 0,1,3\mod 4$.
As a consequence, both the total Milnor invariants $\mu_n\colon \W^\iinfty_n\to \sD_n$  and the realization maps $R^\iinfty_n : \cT^\iinfty_n \to\W^\iinfty_n$ are isomorphisms for these orders.
\end {thm}

The consequences listed in the second statement follow from the fact that $\eta_n$ is the composition
$$
\xymatrix{
\eta_n: \cT^\iinfty_n \ar@{->>}[r]^{R_n^\iinfty} & \W^\iinfty_n \ar@{->>}[r]^{\mu_n} & \sD_{n}
}
$$

Theorem~\ref{thm:LC} is also instrumental in determining the only possible remaining obstructions to computing $\W^\iinfty_{4k-2}$:

%Proposition~6? of Paper 2:
\begin {prop}[\cite{CST4}]\label{prop:kerEta4k-2}
The map sending a rooted tree $J$ to $(J,J)^\iinfty\in\cT^\iinfty_{4k-2}$ induces an isomorphism
\[
\Z_2 \otimes \sL_k \cong\Ker(\eta_{4k-2})
\]
\end {prop}
These symmetric $\iinfty$-trees $(J,J)^\iinfty$ correspond to twisted Whitney disks, and determine
the higher-order Arf invariants $\Arf_k$. All of our above conjectures are equivalent to
the statement that $\W^\iinfty_{4k-2}$ is isomorphic to $\sD_{4k-2}\oplus (\Z_2 \otimes \sL_k)$ via these maps.

Theorem~\ref{thm:Milnor-twisted} and Proposition~\ref{prop:kerEta4k-2} imply Theorem~\ref{thm:mu-epimorphism-on-W-infty} and Corollary~\ref{cor:twisted-classification} above, and \cite{CST4} describes analogous implications of the above-described results in the framed setting (Theorems~\ref{thm:framed-mu-even-and-odd}, \ref{thm:SL}, \ref{thm:framed-Arf}, and Corollary~\ref{cor:mu-sl-arf-classify}).

%%%%%%%%%%%%%%%%%%%%%%%%%%%%%%%%%%%%%%%%%%%

\section{Framed versus twisted Whitney towers}\label{subsec:framed-vs-twisted}
This section describes how the higher-order Sato-Levine and Arf invariants can be interpreted as obstructions to framing a twisted Whitney tower.
The starting point is the following surprisingly simple relation between twisted and framed Whitney towers of various orders:

\begin {prop}[\cite{CST1,CST4}]\label{prop:exact sequence}
For any $n\in\N$, there is a commutative diagram of exact sequences
\[
\xymatrix{ 
0\ar[r] & \W_{2n} \ar[r] &  \W^\iinfty_{2n} \ar[r] & \W_{2n-1}\ar[r] &  \W^\iinfty_{2n-1}\ar[r] & 0 \\
0\ar[r] & \widetilde{\cT}_{2n} \ar@{->>}[u]^{\widetilde{R}_{2n}} \ar[r] &  \cT^\iinfty_{2n}\ar@{->>}[u]^{R^\iinfty_{2n}} \ar[r] & \widetilde\cT_{2n-1}   \ar@{->>}[u]^{\widetilde{R}_{2n-1}} \ar[r] &  \cT^\iinfty_{2n-1}  \ar@{->>}[u]^{R^\iinfty_{2n-1}} \ar[r] & 0
 } \]
 Moreover, there are isomorphisms
\[
\Cok ( \cT_{2n} \to  \cT^\iinfty_{2n}) \cong \Z_2 \otimes \sL'_{n+1} \cong \ker(  \widetilde\cT_{2n-1}  \to  \cT^\iinfty_{2n-1} )
\]
\end {prop}

In the first row, all maps are induced by the identity on the set of links. 
To see the exactness, observe that there is a natural inclusion $\bW_{n} \subseteq \bW_{n}^\iinfty$, and by definition $\bW_{2n-1}^\iinfty = \bW_{2n-1}$.
One then needs to show that indeed 
$\bW_{2n}^\iinfty \subseteq \bW_{2n-1}^\iinfty$, which is accomplished in \cite{CST1}, and then the exact sequence in 
Proposition~\ref{prop:exact sequence} follows since $\W_n := \bW_n/ \bW_{n+1}$ and $\W^\iinfty_n := \bW^\iinfty_n/  \bW^\iinfty_{n+1}$. 

If our above conjectures hold, then for every $n$ the various (vertical) realization maps in the above diagram are isomorphisms, which would lead to a computation of the cokernel and kernel of the map $\W_n\to \W^\iinfty_n$. 
As a consequence, we would obtain new concordance invariants with values in 
$\Z_2 \otimes \sL'_{n+1}$ and defined on $\bW^\iinfty_{2n}$, as the obstructions for a link to bound a framed Whitney tower of order~$2n$. In fact \cite{CST4}, the above-defined higher-order Sato-Levine invariants detect the quotient $\Z_2 \otimes \sL_{n+1}$ of $\Z_2 \otimes \sL'_{n+1}$.
 Levine \cite{L2} showed that the squaring map $X\mapsto [X,X]$ induces an isomorphism
$$
\Z_2 \otimes \sL_{k}
\cong\ker(\Z_2 \otimes \sL'_{2k}\twoheadrightarrow\Z_2 \otimes \sL_{2k}), 
$$
which leads to our proposed higher-order Arf invariants $\Arf_k$.

It is interesting to note that the case $n=0$ leads to the prediction
$
\Cok ( \W_{0} \to  \W^\iinfty_{0}) \cong \Z_2 \otimes \sL_{1} \cong (\Z_2)^m
$
This is indeed the group of framed $m$-component links modulo those with even framings! In fact, the consistency of this computation was the motivating factor to consider filtrations of the set of {\em framed} links $\bL$, rather than just oriented links.

%%%%%%%%%%%%%%%%%%%%%%%%%%%%%%%%%%%%%%%%

\section{Filtrations of homology cylinders}
Let $\Sigma_{g,1}$ denote the compact orientable surface of genus $g$ with one 
boundary component. A \emph{homology cylinder} over $\Sigma_{g,1}$ is a compact $3$ manifold $M$ which is homologically equivalent to the cylinder $\Sigma_{g,1}\times[0,1]$, equipped with standard parameterizations of the two copies of $\Sigma_{g,1}$ at each ``end." 
Two homology cylinders $M_0$ and $M_1$ are said to be \emph{homology cobordant} if there is a compact oriented $4$-manifold $N$ with $\partial N=M_0\cup_{\Sigma_{g,1}} (-M_1)$, such that the inclusions $M_i\hookrightarrow N$ are homology isomorphisms. Let $\mathcal H_g$ be the set of homology cylinders up to homology cobordism over $\Sigma_{g,1}$. $\mathcal H_g$ is a group via the ``stacking" operation. 

Adapting the usual string link definition, Garoufalidis and Levine \cite{GL} introduced an Artin-type representation
$\sigma_n\colon \mathcal{H}_g\to A_0(F/F_{n+1})$ where $F$ is the free group on $2g$ generators, and $A_0(F/F_{n+1})$ is the group of automorphisms $\phi$ of $F/F_{n+1}$ such that  $\phi$ fixes the product $[x_1,y_1]\cdots [x_g,y_g]$ modulo $F_{n+1}$. Here $\{x_i,y_i\}_{i=1}^g$ is a standard symplectic basis for $\Sigma_{g,1}$.  The Johnson (relative weight) filtration of $\mathcal H_g$ is defined by $\mathbb {J}_n=\ker \sigma_n$.  Define the associated graded group ${\sf J}_n=\mathbb J_n/\mathbb J_{n+1}$. Levine showed in \cite{L1} that ${\sf J}_n\cong{\sf D_n}$.

On the other hand, there is a filtration related to Goussarov-Habiro's theory of finite type 3-manifold invariants. One defines the relation of $A_n$-equivalence to be generated by the following move: $M\sim_n M'$ if $M'$ is 
diffeomorphic to $M_C$, for some connected clasper $C$ with $n$ nodes. Let $\mathbb Y_n$ be the subgroup of $\mathcal H_g$ of all homology cylinders $A_n$-equivalent to the trivial one, and let ${\sf Y}_n=\mathbb Y_n/\sim_{n+1}$. Rationally, Levine showed the associated graded groups for these two filtrations to be the same, and are even classified by the tree group $\mathcal T_{n}$:

\begin {thm}[Levine] Levine's map $\eta_{n}$ is the composition\\
\centerline{$\xymatrix{
\cT_{n}\ar@{->>}[r]&{\sf Y}_n\ar[r]&\mathcal {\sf J}_n\ar[r]^\cong &{\sf D}_n
}$} All three maps are rational isomorphisms.
\end {thm}

The story is more subtle over the integers. Using the algebraic methods that we have developed for higher-order intersections we are able to understand how things behave integrally in $3/4$ of the cases, and modulo the question of the non-triviality of higher-order Arf invariants in the remaining cases:

\begin {thm}[\cite{CST5}] \label{thm:filt} For all $k\geq 1$, there are exact sequences
\begin{enumerate}
\item $0\to {\sf Y}_{2k}\to {\sf J}_{2k}\to \mathbb Z_2\otimes{\sf L}_{k+1}\to 0$
\item  $0\to \mathbb Z_2\otimes{\sf L}_{2k+1}\to {\sf Y}_{4k-1}\to{\sf J}_{4k-1}\to 0$
\item $0\to{\sf K}^{\sf Y}_{4k-3}\to{\sf Y}_{4k-3}\to {\sf J}_{4k-3}\to 0$
\item $\mathbb Z_2\otimes {\sf L}_{k}\overset{a_{k}}\to{\sf K}^{\sf Y}_{4k-3}\to \mathbb Z_2\otimes {\sf L}_{2k}\to 0$.
\end{enumerate}
\end {thm}

The calculation of the kernel ${\sf K}^{\sf Y}_{4k-3}$ is thus reduced to the calculation of $\ker(a_{k})$. This is the precise analog of the questions for Whitney towers  about the injectivity of $\alpha_{k}$  and nontriviality of higher-order Arf invariants.

\begin {conj}
The homomorphisms $a_{k}$ are injective for all $k\geq 1$ and hence ${\sf K}^{\sf Y}_{4k-3}\cong \Z_2\otimes {\sf L}'_{2k}$.
\end {conj}

The proof of Theorem~\ref{thm:filt} depends crucially on the resolution of the Levine conjecture, as well as geometric arguments showing that the framing relations in $\widetilde{\cT}$ are also present for homology cylinders.

\vspace{.5em}
{\bf Acknowledgments:}
This paper was written while the first two authors were visiting  Max Planck Institut f\"ur Mathematik in Bonn. They thank MPIM for its stimulating research environment and generous support. The last author was also supported by NSF grants DMS-0806052 and DMS-0757312.

%% PNAS does not support submission of supporting .tex files such as BibTeX.
%% Instead all references must be included in the article .tex document. 
%% If you currently use BibTeX, your bibliography is formed because the 
%% command \verb+\bibliography{}+ brings the <filename>.bbl file into your
%% .tex document. To conform to PNAS requirements, copy the reference listings
%% from your .bbl file and add them to the article .tex file, using the
%% bibliography environment described above.  

%%  Contact pnas@nas.edu if you need assistance with your
%%  bibliography.

% Sample bibliography item in PNAS format:
%% \bibitem{in-text reference} comma-separated author names up to 5,
%% for more than 5 authors use first author last name et al. (year published)
%% article title  {\it Journal Name} volume #: start page-end page.
%% ie,
% \bibitem{Neuhaus} Neuhaus J-M, Sitcher L, Meins F, Jr, Boller T (1991) 
% A short C-terminal sequence is necessary and sufficient for the
% targeting of chitinases to the plant vacuole. 
% {\it Proc Natl Acad Sci USA} 88:10362-10366.

%% Enter the largest bibliography number in the facing curly brackets
%% following \begin{thebibliography}

\end{document}